%
%
%
%
\hsize=5in
\baselineskip=12pt
\vsize=20cm
\parindent=10pt
\pretolerance=40
\predisplaypenalty=0
\displaywidowpenalty=0
\finalhyphendemerits=0
\hfuzz=2pt
\frenchspacing
\footline={\ifnum\pageno=1\else\hfil\tenrm\number\pageno\hfil\fi}
%
%
\input amssym.def
\font\ninerm=cmr9
\font\ninebf=cmbx9
\font\ninei=cmmi9
\skewchar\ninei='177
\font\nineit=cmti9
\def\reffonts{\baselineskip=0.9\baselineskip
	\textfont0=\ninerm
	\def\rm{\fam0\ninerm}%
	\textfont1=\ninei
  \def\bf{\ninebf}%
	\def\it{\nineit}%
	}
%
%
\def\frontmatter{\vbox{}\vskip1cm\bgroup
	\leftskip=0pt plus1fil\rightskip=0pt plus1fil
	\parindent=0pt
	\parfillskip=0pt
	\pretolerance=10000
	}
\def\endfrontmatter{\egroup\bigskip}
\def\title#1{{\bf #1\par}}
\def\author#1{\bigskip#1\par}

\def\thanks#1{\footnote{}{\reffonts\rm\noindent#1\hfil}}
\def\section#1\par{\ifdim\lastskip<\bigskipamount\removelastskip\fi
	\penalty-250\bigskip
	\vbox{\leftskip=0pt plus1fil\rightskip=0pt plus1fil
	\parindent=0pt
	\parfillskip=0pt
  \pretolerance=10000{\bf#1}}\nobreak\medskip
	}
\def\emph#1{{\it #1}\/}
\def\proclaim#1. {\medbreak\bgroup{\noindent\bf#1.}\ \it}
\def\endproclaim{\egroup
	\ifdim\lastskip<\medskipamount\removelastskip\medskip\fi}
\newdimen\itemsize
\def\setitemsize#1 {{\setbox0\hbox{#1\ }
	\global\itemsize=\wd0}}
\def\item#1 #2\par{\ifdim\lastskip<\smallskipamount\removelastskip\smallskip\fi
	{\leftskip=\itemsize
	\noindent\hskip-\leftskip
	\hbox to\leftskip{\hfil\rm#1\ }#2\par}\smallskip}
\def\Proof#1. {\ifdim\lastskip<\medskipamount\removelastskip\medskip\fi
	{\noindent\it Proof\if\space#1\space\else\ \fi#1.}\ }
\def\endproof{\hfill\hbox{}\quad\hbox{}\hfill\llap{$\square$}\medskip}
%
%
\newcount\citation
\newtoks\citetoks
\def\citedef#1\endcitedef{\citetoks={#1\endcitedef}}
\def\endcitedef#1\endcitedef{}
\def\citenum#1{\citation=0\def\curcite{#1}%
	\expandafter\checkendcite\the\citetoks}
\def\checkendcite#1{\ifx\endcitedef#1?\else
	\expandafter\lookcite\expandafter#1\fi}
\def\lookcite#1 {\advance\citation by1\def\auxcite{#1}%
	\ifx\auxcite\curcite\the\citation\expandafter\endcitedef\else
	\expandafter\checkendcite\fi}
\def\cite#1{\makecite#1,\cite}
\def\makecite#1,#2{[\citenum{#1}\ifx\cite#2]\else\expandafter\clearcite\expandafter#2\fi}
\def\clearcite#1,\cite{, #1]}
%
%
\def\references{\section References\par
	\bgroup
	\parindent=0pt
	\reffonts
	\rm
	\frenchspacing
	\setbox0\hbox{99. }\leftskip=\wd0
	}
\def\endreferences{\egroup}
\newtoks\authtoks
\newif\iffirstauth
\def\checkendauth#1{\ifx\auth#1%
    \iffirstauth\the\authtoks
    \else{} and \the\authtoks\fi,%
  \else\iffirstauth\the\authtoks\firstauthfalse
    \else, \the\authtoks\fi
    \expandafter\nextauth\expandafter#1\fi
	}
\def\nextauth#1,#2;{\authtoks={#1 #2}\checkendauth}
\def\auth#1{\nextauth#1;\auth}
\newif\ifbookref
\def\nextref#1 {\par\hskip-\leftskip
	\hbox to\leftskip{\hfil\citenum{#1}.\ }%
	\initnextref}
\def\initnextref{\bookreffalse\firstauthtrue\ignorespaces}
\def\paper#1{{\it#1},}
\def\book#1{\bookreftrue{\it#1},}
\def\journal#1{#1}
\def\Vol#1{{\bf#1}}
\def\publisher#1{#1,}
\def\Year#1{\ifbookref #1.\else(#1)\fi}
\def\Pages#1{\makepages#1.}
\long\def\makepages#1-#2.#3{#1--#2\ifx\par#3.\fi#3}
\def\inRus{{ \rm(in Russian)}}
\def\etransl#1{English translation in \journal{#1}}
%
%
\newsymbol\square 1003
\newsymbol\bbk 207C
\newsymbol\barwedge 125A
\let\wdg\barwedge
\let\ot\otimes
\let\sbs\subset
\let\<\langle
\let\>\rangle
\def\defop#1#2{\def#1{\mathop{\rm #2}\nolimits}}
\defop\Alt{Alt}
\defop\chr{char}
\defop\End{End}
\defop\gl{\frak{gl}}
\defop\GL{GL}
\defop\HeckeSym{HeckeSym}
\defop\Id{Id}
\defop\Im{Im}
\defop\Ker{Ker}
\def\mod#1{\mskip12mu(\mathop{\rm mod}#1)}
\defop\rk{rank}
\defop\sl{\frak{sl}}
\defop\tr{tr}
\let\al\alpha
\let\la\lambda
\let\om\omega
\def\omt{\widetilde\om}
\let\ph\varphi
\let\si\sigma
\let\th\theta
\let\De\Delta
\def\bbS{{\Bbb S}}
\def\bbT{{\Bbb T}}

\citedef
Ch-P
Dr83
Ew-O94
Ger-Gi98
Gur90
Sto91
\endcitedef

\frontmatter
\title{Hecke symmetries associated with the polynomial algebra\break 
in 3 commuting indeterminates}
\author{Serge Skryabin}
\endfrontmatter

All Hecke symmetries corresponding to quantum groups of the $\GL(3)$ family can 
be found by solving a large set of equations in the coordinates of two cubic 
tensors associated with a pair of Artin-Schelter regular graded algebras of 
global dimension 3. This method worked out by Ewen and Ogievetsky \cite{Ew-O94} 
requires a case by case analysis. Each case is distinguished by the choice of 
one cubic tensor and a linear operator whose action should fix that tensor. 

The alternating cubic tensors correspond to the free commutative algebras with 
3 generators. This case appears to be somewhat troublesome in the approach of 
Ewen and Ogievetsky as it involves no restriction on the linear operators 
which extend to automorphisms of the algebra. In fact, the necessary equations 
with this choice of a cubic tensor were not investigated in \cite{Ew-O94}. 
Instead it was proposed to interchange the roles of two tensors and look at 
all solutions where the second tensor was found to be alternating.

In the present paper we pursue quite a different approach which enables us to 
derive a general formula for all Hecke symmetries $R$ with the $R$-symmetric 
algebra $\bbS(V,R)$ being the polynomial algebra in 3 commuting indeterminates. 
It will be shown that such a Hecke symmetry is determined by a bivector and a 
symmetric bilinear form on a 3-dimensional vector space $V$. Precise statements 
of results are given in section 5 of the paper.

Rather than trying to fix all possibilities for the second cubic tensor we 
reformulate the braid equation for $R$ in terms of a certain collection of 
linear functions $\ell_{xy}\in V^*$ indexed by pairs of vectors $x,y\in V$. In 
this form equations are easily manageable. Especially, it will be immediately 
clear that all $\ell_{xx}$, $x\in V$, are scalar multiples of one function. 
This property will further lead to the determination of $\ell_{xy}$ for 
arbitrary $x,y$ and the recovery of the Hecke symmetry itself. At the end we 
will still need some checks to ensure that the braid equation holds in full 
generality. 

Our solution exploits invariance of the ordinary symmetric algebra $\bbS(V)$ 
and some peculiarities of the 3-dimensional space. Unfortunately, this method 
does not generalize to higher dimensions.

As a byproduct we will see that each Hecke symmetry with the prescribed 
algebra $\bbS(V,R)$ is a deformation of the flip operator $R_0$ that sends 
$x\ot y$ to $y\ot x$ for all $x,y\in V$. This was by no means clear a priori. 
The already mentioned preprint \cite{Ew-O94} gives examples of $R$-matrices 
not obtained by such a deformation.

In section 6 we determine the equivalence classes of Hecke symmetries with the 
prescribed algebra $\bbS(V,R)$. There are two families with an arbitrary 
parameter $q\ne\nobreak1$ of the quadratic Hecke relation and six equivalence 
classes with $q=1$. The latter can be obtained as specializations of the 
operators in the 1-parameter families, but we list all cases with $q=1$ as 
separate types. We provide explicit formulas for the Hecke symmetry $R$ in 
each type and for the corresponding classical $r$-matrix.

As a matter of comparison, Ewen and Ogievetsky mention a solution with an 
alternating cubic tensor only once \cite{Ew-O94, p. 18}. This does not tell the 
full story as the twists used there distort the tensor and the corresponding 
algebra. One can search for the other solutions by twisting back.

\section
1. Notation

Let $V$ be a vector space of dimension 3 over a field $\bbk$. In the tensor 
algebra $\bbT(V)$ we have
$$
xy=x\ot y,\qquad xyz=x\ot y\ot z
$$
for $x,y,z\in V$. This will be understood in all formulas. For each $k\ge0$ 
denote by $\Alt_k$ the subspace of alternating tensors in $V^{\ot k}$. There 
is a multiplication $\wdg$ which makes the direct sum of spaces $\Alt_k$ an 
associative algebra isomorphic to the exterior algebra $\bigwedge V$. Explicitly,
$$
x\wdg y=xy-yx,\qquad x\wdg y\wdg z=xyz+yzx+zxy-zyx-xzy-yxz
$$
for $x,y,z\in V$. Essentially, we may identify $x\wdg y$ and $x\wdg y\wdg z$ 
with the bivector $x\wedge y$ and the trivector $x\wedge y\wedge z$.

Fix a nonzero alternating trilinear form $\om:V\times V\times V\to\bbk$. Since 
the space $\Alt_3$ has dimension 1, there is a linear bijection 
$\omt:\Alt_3\to\bbk$ such that
$$
\omt(x\wdg y\wdg z)=\om(x,y,z)\quad\hbox{for $ x,y,z\in V$}.
$$
Let $\om_{xy}\in V^*$ be the linear forms defined by the rule
$$
\om_{xy}(v)=\om(x,y,v),\qquad x,y,v\in V.
$$
Note that $\om_{xy}\ne0$ if and only if $x$ and $y$ are linearly independent, 
and in this case $\Ker\om_{xy}$ is spanned by $x$ and $y$. If $l\in V^*$ is 
any nonzero linear form, then $l=\om_{ab}$ for a suitably normalized basis 
$a,b$ of the 2-dimensional subspace $\Ker l\sbs V$. Thus 
$V^*=\{\om_{xy}\mid x,y\in V\}$.

\section
2. Reformulation of the braid equation

Let $0\ne q\in\bbk$. According to Gurevich \cite{Gur90} a \emph{Hecke symmetry} 
with parameter $q$ on $V$ is any linear operator $R:V\ot V\to V\ot V$ 
satisfying the braid equation
$$ 
(R\ot\Id_V)(\Id_V\ot\,R)(R\ot\Id_V)=(\Id_V\ot\,R)(R\ot\Id_V)(\Id_V\ot\,R) 
$$
and the quadratic Hecke relation
$$
(R-q\cdot\Id_{V\ot V})(R+\Id_{V\ot V})=0.
$$
The \emph{$R$-symmetric algebra} $\bbS(V,R)$ is the factor algebra of the 
tensor algebra $\bbT(V)$ by the graded ideal generated by the subspace 
$\,\Im Y\sbs V^{\ot 2}\,$ where
$$
Y=q\cdot\Id_{V\ot V}-R\eqno(2.1)
$$
is the \emph{$R$-skewsymmetrizer}.

Our aim is to describe all Hecke symmetries $R\in\GL(V^{\ot2})$ such that 
$\bbS(V,R)$ is the ordinary symmetric algebra of $V$, i.e., 
$\bbS(V,R)=\bbT(V)/I$ where $I$ is the ideal of $\bbT(V)$ generated by the 
3-dimensional subspace $\Alt_2$ of alternating tensors in $V^{\ot2}$.

So let $R$ be such a Hecke symmetry with parameter $q$. Then $\,\Im Y=\Alt_2$. 
The Hecke relation for $R$ implies that $Y^2=(q+1)Y$. Hence
$$
Yw=(q+1)w\qquad\hbox{for all $w\in\Alt_2$}.\eqno(2.2)
$$
The braid equation for $R$ is equivalent to the following equation for $Y$:
$$
\eqalign{
&(\Id_V\ot\,Y)(Y\ot\Id_V)(\Id_V\ot\,Y)-q\,(\Id_V\ot\,Y)\cr
&\hskip30mm=(Y\ot\Id_V)(\Id_V\ot\,Y)(Y\ot\Id_V)-q\,(Y\ot\Id_V).
}\eqno(2.3)
$$
Here the linear operators $\Id\ot Y$ and $Y\ot\Id$ acting on $V^{\ot3}$ have 
images, respectively, $V\ot\Alt_2$ and $\Alt_2\ot V$. Therefore the two equal 
operators in the above equality have images in the 1-dimensional subspace
$$
\Alt_3=(V\ot\Alt_2)\cap(\Alt_2\ot V)
$$
of alternating tensors in $V^{\ot3}$. So it follows that
$$
(\Id_V\ot\,Y)(Y\ot\Id_V)w-qw\in\Alt_3\qquad\hbox{for all $w\in V\ot\Alt_2$}.\eqno(2.4)
$$

We will eventually find a general form of such operators $Y$. First we are 
going to reformulate the requested condition in a different form starting with 
its coordinate representation as an intermediate step. Let $e_1,e_2,e_3$ be a 
linear basis for $V$, and let $Y_{ij}^{kl}\in\bbk$ be the coefficients in the 
expressions
$$
Y(e_ie_j)=\sum Y_{ij}^{kl}e_ke_l,\qquad 
Y_{ij}^{kl}=-Y_{ij}^{lk},\quad Y_{ij}^{kk}=0.
$$
Here and later we assume that the summation is over the indices repeated as 
subscripts and superscripts. Then
$$
(\Id\ot\,Y)(Y\ot\Id)(e_ie_je_k-e_ie_ke_j)
=\sum(Y_{ij}^{rl}Y_{lk}^{st}-Y_{ik}^{rl}Y_{lj}^{st})e_re_se_t.
$$
Since $e_ie_je_k-e_ie_ke_j\in V\ot\Alt_2$, we must have
$$
\sum(Y_{ij}^{rl}Y_{lk}^{st}-Y_{ik}^{rl}Y_{lj}^{st})e_re_se_t
\equiv q(e_ie_je_k-e_ie_ke_j)\mod{\Alt_3}.
$$
If $s=r$, then the basis monomial $e_re_se_t$ has zero coefficient in 
the elements of $\Alt_3$. Comparison of coefficients yields
$$
\sum\,\bigl(Y_{ij}^{rl}Y_{lk}^{rt}-Y_{ik}^{rl}Y_{lj}^{rt}\bigr)=\cases{
\hphantom{-}0&if $i\ne r$, or $t=r$, or $\{j,k\}\ne\{r,t\}$,\cr
\hphantom{-}q&if $i=j=r$ and $k=t\ne r$,\cr
-q&if $i=k=r$ and $j=t\ne r$
}\eqno(2.5)
$$
where the summation is over $l$. We will not need to look at the basis 
monomials $e_re_se_t$ with $s\ne r$ due to the following observation:

\proclaim
Lemma 2.1.
Let $Y:V^{\ot2}\to V^{\ot2}$ be a linear operator such that\/ $\Im Y\sbs\Alt_2$. 
For the containments $(2.4)$ to hold it is necessary and sufficient that in every 
basis of $V$ the matrix components of $Y$ satisfy $(2.5)$.
\endproclaim

\Proof.
Equalities (2.5) mean that the image of the linear operator
$$
Z=\bigl((\Id_V\ot Y)(Y\ot\Id_V)-q\cdot\Id_{V^{\ot3}}\bigr)\Big|_{V\ot\Alt_2}
$$
is contained in the 6-dimensional linear subspace $W_{e_1,e_2,e_3}$ of 
$V^{\ot3}$ spanned by the basis monomials $e_re_se_t$ with $r,s,t$ pairwise 
distinct. If these equalities hold in every basis of $V$, then $\Im Z\sbs W$ 
where $W$ is the intersection of those subspaces $W_{e_1,e_2,e_3}$ taken for 
all bases of $V$. Since $W$ is a $\GL(V)$-submodule of $V^{\ot3}$, it is easy 
to see that $W=\Alt_3$. Thus $\,\Im Z\sbs\Alt_3$.
\endproof

Next we will find a coordinate-free interpretation of equalities (2.5). 
The multiplication $\wdg$ produces the elements $x\wdg Y(yz)\in\Alt_3$ for 
$x,y,z\in V$. So there are linear forms $\ell_{xy}\in V^*$ defined by the rule
$$
\ell_{xy}(z)=\omt\bigl(x\wdg Y(yz)\bigr),\qquad x,y,z\in V.\eqno(2.6)
$$
Note that the expression $\ell_{xy}(z)$ is a trilinear function of $(x,y,z)$. 
For $f,g\in V^*$ we identify $f\wedge g\in\bigwedge^2V^*$ with an alternating 
bilinear form on $V$ setting
$$
(f\wedge g)(u,v)=f(u)g(v)-f(v)g(u),\qquad u,v\in V.
$$

\proclaim
Lemma 2.2.
Let $Y:V^{\ot2}\to V^{\ot2}$ be a linear operator such that $\Im Y\sbs\Alt_2$. 
For the containments $(2.4)$ to hold it is necessary and sufficient that
$$
(\ell_{xy}\wedge\ell_{xz}-\ell_{xx}\wedge\ell_{yz})(u,v)=q\,\om(x,y,z)\,\om(x,u,v)\eqno(2.7)
$$
for all $x,y,z,u,v\in V$.
\endproclaim

\Proof.
The expressions in the left and right hand sides of (2.7) are linear in $y,z,u,v$. 
If either $x=0$ or $y=x$, then both sides are equal to $0$. So we may assume 
that $x$ and $y$ are linearly independent. Let $e_1,e_2,e_3$ be any basis of $V$ 
such that $x=e_2$ and $y=e_3$. Put $\al=\om(e_1,e_2,e_3)$. So $0\ne\al\in\bbk$. 
By linearity it suffices to check (2.7) for all triples $z,u,v$ consisting of 
basis elements. Let $z=e_i$. Since
$$
e_2\wdg Y(e_ie_j)=Y_{ij}^{13}e_2\wdg e_1\wdg e_3,\qquad
e_3\wdg Y(e_ie_j)=Y_{ij}^{12}e_1\wdg e_2\wdg e_3,      
$$
we get
$$
\ell_{xz}(e_j)=\omt\bigl(e_2\wdg Y(e_ie_j)\bigr)=-\al Y_{ij}^{13},\qquad
\ell_{yz}(e_j)=\omt\bigl(e_3\wdg Y(e_ie_j)\bigr)=\al Y_{ij}^{12}
$$
for each $j=1,2,3$. With $i=2$ and $i=3$ the first of the above 
equalities gives
$$
\ell_{xx}(e_j)=-\al Y_{2j}^{13},\qquad\ell_{xy}(e_j)=-\al Y_{3j}^{13}.
$$
Hence
$$
\eqalign{
(\ell_{xy}\wedge\ell_{xz}-\ell_{xx}\wedge\ell_{yz})(e_k,e_j)
&=\al^2(Y_{ij}^{13}Y_{3k}^{13}-Y_{ik}^{13}Y_{3j}^{13})
+\al^2(Y_{ij}^{12}Y_{2k}^{13}-Y_{ik}^{12}Y_{2j}^{13})\cr
&=\al^2\sum\,\bigl(Y_{ij}^{1l}Y_{lk}^{13}-Y_{ik}^{1l}Y_{lj}^{13}\bigr)
}
$$
(the $l=1$ term in the sum vanishes since $Y_{ij}^{11}=Y_{ik}^{11}=0$). 
On the other hand,
$$
\om(e_2,e_3,e_i)\,\om(e_2,e_k,e_j)=\cases{
\hphantom{-}0&if either $i\ne 1$ or $\{j,k\}\ne\{1,3\}$,\cr
\hphantom{-}\al^2&if $i=j=1$ and $k=3$,\cr
-\al^2&if $i=k=1$ and $j=3$.
}
$$
We see that (2.7) holds for all $x,y,z,u,v\in V$ if and only if (2.5) holds in 
every basis of $V$ for $r=1$, $t=3$ and arbitrary $i,j,k$. For any other pair 
of indices $r\ne t$ equalities (2.5) reduce to the case $r=1$, $t=3$ by 
renumbering the basis elements. If $r=t$, then the sum in the left hand side 
of (2.5) vanishes since $Y_{lj}^{rr}=0$ for all $l,j$. Thus Lemma 2.2 is just a 
reformulation of Lemma 2.1.
\endproof

\section
3. Solution

Now comes the main step. Assume that $\,\chr\bbk\ne2$.

\proclaim
Proposition 3.1.
Let $\,Y:V^{\ot2}\to V^{\ot2}\,$ be a linear operator such that $\,\Im Y\sbs\Alt_2$. 
Then $Y$ satisfies conditions $(2.2)$ and $(2.4)$ if and only if there exist two 
vectors $a,b\in V$ and a symmetric bilinear form $g:V\times V\to\bbk$ such that
$$
(q-1)^2=-4\De\quad\hbox{where $\De=g(a,a)\,g(b,b)-g(a,b)^2$}\eqno(3.1)
$$
and $Y$ is given by the formula
$$
Y(xy)=g(x,y)\,a\wdg b+x\wdg Ty+y\wdg Tx+\smash{q+1\over2}x\wdg y,\qquad x,y\in V,
\eqno(3.2)
$$
where $T:V\to V$ is the linear operator defined by the rule
$$
Tv=g(b,v)\,a-g(a,v)\,b,\qquad v\in V.\eqno(3.3)
$$
Moreover, such an operator $Y$ also satisfies $(2.3)$.
\endproclaim

\Proof.
Since the assignment $(x,t)\mapsto\omt(x\wdg t)$ defines a nondegenerate 
bilinear pairing between the vector spaces $V$ and $\Alt_2$, it is seen 
from (2.6) that the operator $Y$ is uniquely determined by the collection of 
linear forms $\ell_{xy}$. We have
$$
\omt\bigl(x\wdg Y(yz-zy)\bigr)=\ell_{xy}(z)-\ell_{xz}(y),\qquad x,y,z\in V.
$$
Since the space $\Alt_2$ is spanned by the tensors $y\wdg z=yz-zy$, it 
follows that condition (2.2) amounts to the identity
$$
\ell_{xy}(z)-\ell_{xz}(y)=(q+1)\,\om(x,y,z).\eqno(3.4)
$$
By Lemma 2.2 condition (2.4) is expressed by means of identity (2.7). We will 
determine the forms $\ell_{xy}$ satisfying (2.7) and (3.4), and then derive a 
formula for $Y$.

We first use a simple consequence of (2.7). Substituting $z=y$ in (2.7) we get
$$
\ell_{xx}\wedge\ell_{yy}=0\quad\hbox{for all $x,y\in V$}
$$
since $\ell_{xy}\wedge\ell_{xy}=0$ and $\om(x,y,y)=0$. This means that all nonzero 
functions in the set $\{\ell_{xx}\mid x\in V\}$ are scalar multiples of each other. 
Note also that the assignment $x\mapsto\ell_{xx}$ gives a homogeneous polynomial 
map of degree 2. Hence there is a pair $l$, $\ph$ consisting of a linear form 
and a quadratic form on $V$ such that
$$
\ell_{xx}=\ph(x)\,l=g(x,x)\,l\qquad\hbox{for all $x\in V$}\eqno(3.5)
$$
where $g$ is the unique symmetric bilinear form such that $\ph(x)=g(x,x)$ for 
all $x$. If $\ell_{xx}=0$ for all $x$, then we may use $l=0$ and $g=0$. Let 
$a,b\in V$ be such that $l=\om_{ab}$. So
$$
l(v)=\om(a,b,v)\qquad\hbox{for all $v\in V$}.\eqno(3.6)
$$

Formula (3.4) with $y=x$ gives $\ell_{xz}(x)=\ell_{xx}(z)$. Hence
$$
\ell_{yx}(y)=\ell_{yy}(x)=g(y,y)\,l(x).
$$
Since $\,\ell_{xy}+\ell_{yx}=2g(x,y)l\,$ by linearization of (3.5), we get
$$
\ell_{xy}(y)=2g(x,y)\,l(y)-\ell_{yx}(y)=2g(x,y)\,l(y)-g(y,y)\,l(x)
$$
and
$$
\ell_{xy}(z)+\ell_{xz}(y)=2g(x,y)\,l(z)+2g(x,z)\,l(y)-2g(y,z)\,l(x).
$$
by linearization. Combined with (3.4) the last identity yields
$$
\ell_{xy}(z)=g(x,y)\,l(z)+g(x,z)\,l(y)-g(y,z)\,l(x)+{q+1\over2}\om(x,y,z).\eqno(3.7)
$$
It is readily seen that for any $l$ and $g$ the linear forms $\ell_{xy}$ 
defined by (3.7) satisfy (3.4). We will show that (2.7) is satisfied if and only 
if (3.1) holds. This step is somewhat harder.

Since $\bigwedge^4V=0$, every alternating multilinear form of 4 arguments in 
$V$ vanishes. It follows that for any $\xi\in V^*$ and $v_1,v_2,v_3,v_4\in V$ 
there is an equality
$$
\eqalign{
\xi(v_1)\,\om(v_2,v_3,v_4)&-\xi(v_2)\,\om(v_1,v_3,v_4)\cr
&+\xi(v_3)\,\om(v_1,v_2,v_4)-\xi(v_4)\,\om(v_1,v_2,v_3)=0.
}\eqno(3.8)
$$
Taking $\xi$ such that $\xi(v)=g(v,z)$ and using (3.6), we get
$$
\eqalign{
g(x,z)\,l(y)-g(y,z)\,l(x)&=g(b,z)\,\om(a,x,y)-g(a,z)\,\om(b,x,y)\cr
&=\om(Tz,x,y)=\om(x,y,Tz).
}\eqno(3.9)
$$
Thus
$$
\ell_{xy}(z)=g(x,y)\,l(z)+\om(x,y,Tz)+{q+1\over2}\om(x,y,z)\eqno(3.10)
$$
for all $x,y,z$, i.e.,
$$
\ell_{xy}=g(x,y)l+\om_{xy}T+{q+1\over2}\om_{xy}
=g(x,y)l+\om_{xy}T'
$$
where we put $\,T'=T+\displaystyle{q+1\over2}\Id\,$. Now
$$
\ell_{xy}\wedge\ell_{xz}-\ell_{xx}\wedge\ell_{yz}=l\wedge hT' 
+\om_{xy}T'\wedge\om_{xz}T'
$$
where $\,h=g(x,y)\,\om_{xz}-g(x,z)\,\om_{xy}-g(x,x)\,\om_{yz}\,$. As a special 
case of (3.8) we have
$$
h(u)=-g(x,u)\,\om(x,y,z)\quad\hbox{for all $u\in V$}.
$$
Hence
$$
(l\wedge hT')(u,v)
=\bigl(g(x,T'u)\,l(v)-g(x,T'v)\,l(u)\bigr)\,\om(x,y,z).
$$
As another special case of (3.8) we take $\xi=\om_{xy}$ and deduce that
$$
\xi(u)\,\om(x,z,v)-\xi(v)\,\om(x,z,u)=\xi(z)\,\om(x,u,v)
$$
since $\xi(x)=0$. In other words,
$$
(\om_{xy}\wedge\om_{xz})(u,v)=\om(x,y,z)\,\om(x,u,v).
$$
It follows that
$$
(\om_{xy}T'\wedge\om_{xz}T')(u,v)=(\om_{xy}\wedge\om_{xz})(T'u,T'v)=
\om(x,y,z)\,\om(x,T'u,T'v).
$$
Thus the verification of (2.7) reduces to the identity
$$
g(x,T'u)\,l(v)-g(x,T'v)\,l(u)+\om(x,T'u,T'v)=q\,\om(x,u,v).\eqno(3.11)
$$

Let us rewrite the left hand side of (3.11) applying several known identities. 
By (3.9)
$$
g(x,u)\,l(v)-g(x,v)\,l(u)=\om(Tx,u,v).
$$
Since $\,g(x,Tu)=g(u,b)\,g(x,a)-g(u,a)\,g(x,b)=-g(Tx,u)\,$, we also have
$$
g(x,Tu)\,l(v)-g(x,Tv)\,l(u)=-\om(T^2x,u,v).
$$
Hence
$$
\eqalign{
&g(x,T'u)\,l(v)-g(x,T'v)\,l(u)+\om(x,T'u,T'v)\cr
&\qquad={q+1\over2}\bigl(\om(Tx,u,v)+\om(x,Tu,v)+\om(x,u,Tv)\bigr)+{(q+1)^2\over4}\om(x,u,v)\cr
&\qquad\qquad\qquad-\om(T^2x,u,v)+\om(x,Tu,Tv)
}
$$
Since $\,\tr T=0\,$, the first 3 terms sum up to 0, and also
$$
\eqalign{
-\om(T^2x,u,v)+\om(x,Tu,Tv)&=\om(Tx,u,Tv)+\om(Tx,Tu,v)+\om(x,Tu,Tv)\cr
&=c_2(T)\,\om(x,u,v)
}
$$
where $c_2(T)$ is the second coefficient of the characteristic polynomial of $T$.

If $a\wedge b=0$, then $T=0$ and $c_2(T)=0$. Otherwise the image of $T$ 
is contained in the 2-dimensional subspace $\<a,b\>$ of $V$ spanned by $a$ and 
$b$. The restriction of $T$ to this subspace has the matrix
$$
\pmatrix{g(b,a)&g(b,b)\cr -g(a,a)&-g(a,b)},
$$
and so $\,c_2(T)=\det T|_{\<a,b\>}=g(a,a)g(b,b)-g(a,b)^2$. We see that 
$c_2(T)=\De$ in any case. As a consequence, identity (3.11), and therefore also 
(2.7), are satisfied if and only if
$$
{(q+1)^2\over4}+\De=q,
$$
which can be rewritten as (3.1).

Finally we come to the operator $Y$. To derive (3.2) we use (3.9) to obtain 
the expression
$$
g(x,y)\,l(z)=g(y,z)\,l(x)+\om(Ty,x,z).
$$
Now (3.10) can be rewritten as
$$
\displaylines{ 
\ell_{xy}(z)=g(y,z)\,l(x)+\om(x,z,Ty)+\om(x,y,Tz)+{q+1\over2}\om(x,y,z)
=\omt\bigl(x\wdg w)\cr
\hbox{where}\qquad w=g(y,z)\,a\wdg b+y\wdg Tz+z\wdg Ty+{q+1\over2}y\wdg z.
}
$$
Hence $\,Y(yz)=w\,$ as $Y(yz)$ is characterized as the unique element 
$w\in\Alt_2$ such that $\omt\bigl(x\wdg w\bigr)=\ell_{xy}(z)$ for all $x\in 
V$.

It remains to establish (2.3). This will be done separately in the next 
section.
\endproof

\section
4. Verification of the braid equation

We will check that any linear operator $Y:V^{\ot2}\to V^{\ot2}$ given by 
formula (3.2) satisfies (2.3). This equation for $Y$ is equivalent to the braid 
equation for the operator $R=q\cdot\Id-Y$. It would be quite tiresome to 
compute separately the left and right hand sides of (2.3). We will find an 
expression for the difference between the two sides making use of additional 
symmetry properties.

Put $Y_1=Y\ot\Id_V$ and $Y_2=\Id_V\ot Y$ for short. So we have to check that
$$
Y_1Y_2Y_1w-q{\mskip1mu}Y_1w=Y_2Y_1Y_2w-q{\mskip1mu}Y_2w\eqno(4.1)
$$
for all $w\in V^{\ot3}$.

By the construction in the previous section $Y$ satisfies (2.2) and (2.4). The 
operator $Y_1$ gives by restriction a linear map $V\ot\Alt_2\to\Alt_2\ot V$, 
and $Y_2$ gives a linear map in the opposite direction. Both $Y_1$ and $Y_2$ 
act on alternating tensors in $V^{\ot3}$ as the multiplication by $q+1$. Hence 
there are induced linear maps
$$
(V\ot\Alt_2)/\Alt_3\to(\Alt_2\ot V)/\Alt_3\to(V\ot\Alt_2)/\Alt_3
$$
whose composition is $q\cdot\Id$ in view of (2.4). Since the two vector spaces 
involved here have equal dimension, those maps are invertible. Hence the 
composition of the two maps in reversed order is also $q\cdot\Id$, i.e.,
$$
Y_1Y_2w-qw\in\Alt_3\qquad\hbox{for all $w\in\Alt_2\ot V$}.\eqno(4.2)
$$

If $w\in V\ot\Alt_2$, then $Y_2w=(q+1)w$ by (2.2) and $Y_2Y_1w-qw\in\Alt_3$ by 
(2.4). In this case both sides of (4.1) are equal to $(q+1)(Y_2Y_1w-qw)$.
Similarly, both sides of (4.1) are equal to $(q+1)(Y_1Y_2w-qw)$ when 
$w\in\Alt_2\ot V$.

Thus (4.1) holds for all $w\in I_3$ where $I_3=\Alt_2\ot V+V\ot\Alt_2$ is the 
degree 3 homogeneous component of the ideal $I\sbs\bbT(V)$ defining the 
symmetric algebra $\bbS(V)$. So full generality of (4.1) will follow once we 
check this equality for some set of tensors $w\in V^{\ot3}$ whose images in 
$\bbS(V)$ span the vector space $\bbS_3(V)=V^{\ot3}/I_3$.

If $\chr\bbk\ne2,3$ then $\bbS_3(V)$ is spanned by the monomials $x^3$ with 
$x\in V$. This is no longer true when $\chr\bbk=3$, but we note that it 
suffices to check (4.1) in characteristic 0 only. Indeed, taking a commutative 
local domain $K$ with residue field $\bbk$ and the field of fractions $Q(K)$ 
of characteristic 0, one can lift the vectors $a,b$ and the bilinear form $g$ 
to the rank 3 free $K$-module $K^3$. The same formula (3.2) defines then a 
$K$-linear endomorphism of $K^3\ot_K\!K^3$, and the verification that it 
satisfies (4.1) can be done over $Q(K)$.

The proof of (4.1) is thus reduced to the case where $w$ is the tensor 
$x^3\in V^{\ot3}$ for some $x\in V$, i.e., we have to show that
$$
Y_1Y_2(tx)-q{\mskip1mu}tx=Y_2Y_1(xt)-q{\mskip1mu}xt
\quad\hbox{when $t=Y(x^2)\in\Alt_2$}.\eqno(4.3)
$$
Let $\si:V^{\ot3}\to V^{\ot3}$ be the linear operator defined by the formula
$$
\si(xyz)=yzx,\qquad x,y,z\in V.
$$
It maps $V\ot\Alt_2$ onto $\Alt_2\ot V$ and acts as the identity operator on 
$\Alt_3$. Since both sides of (4.3) lie in $\Alt_3$ by (2.4) and (4.2), equality 
(4.3) can be rewritten as
$$
Y_1Y_2(tx)-q\,tx=\si\bigl(Y_2Y_1(xt)-q\,xt\bigr).
$$
Since $\si(xt)=tx$, it is equivalent to the equality
$$
Y_1Y_2(tx)=\si Y_2Y_1(xt).\eqno(4.4)
$$
We will compare the left and right hand sides above by means of the following

\proclaim
Lemma 4.1.
\quad$\,Y_1Y_2(tx)-\si Y_2Y_1(xt)=2\,Tx\wdg t\,$ for all $x\in V\!$ and 
$t\in\Alt_2\,$.
\endproclaim

\Proof.
We may assume that $t=y\wdg z=yz-zy$ for some $y,z\in V$. Then
$$
\eqalign{
Y_1(xt)+\si Y_2(tx)&=Y(xy)z-Y(xz)y+\si\bigl(yY(zx)-zY(yx)\bigr)\cr
&=Y(xy-yx)z-Y(xz-zx)y\cr
&=(q+1)(xyz-yxz-xzy+zxy)=(q+1)(x\wdg t-tx).
}
$$
Here $x\wdg t=x\wdg y\wdg z\in\Alt_3$. Let $\al:\Alt_2\ot V\to\Alt_3$ be the 
linear map such that $\,\al(tx)=t\wdg x=x\wdg t\,$ for $x\in V\!$ and 
$t\in\Alt_2$. We get the identity
$$
Y_1(\si^{-1}w)=-\si Y_2w+(q+1)(\al w-w),\qquad w\in\Alt_2\ot V,
$$
as for $w=tx$ it amounts to the equality written previously. For 
$w=Y_1(xt)$ it gives
$$
Y_1\bigl(\si^{-1}Y_1(xt)\bigr)=-\si Y_2Y_1(xt)
+(q+1)\bigl(Y(xy)\wdg z-Y(xz)\wdg y-Y_1(xt)\bigr).
$$
On the other hand, applying $\si^{-1}$, we find
$$
Y_2(tx)=-\si^{-1}Y_1(xt)+(q+1)(x\wdg t-xt),
$$
and it follows that
$$
\eqalign{
Y_1Y_2(tx)&=-Y_1\bigl(\si^{-1}Y_1(xt)\bigr)+(q+1)^2x\wdg t-(q+1)Y_1(xt)\cr
&=\si Y_2Y_1(xt)-(q+1)Y(xy)\wdg z+(q+1)Y(xz)\wdg y+(q+1)^2x\wdg t.
}
$$
Plugging in the explicit expression (3.2), this gives
$$
\eqalign{
Y_1Y_2(tx)-\si Y_2Y_1(xt)
&=(q+1)\bigl(\,-g(x,y)\,a\wdg b\wdg z+g(x,z)\,a\wdg b\wdg y\cr
&\qquad-x\wdg Ty\wdg z+x\wdg Tz\wdg y-y\wdg Tx\wdg z+z\wdg Tx\wdg y\bigr).
}
$$
Note that\vadjust{\vskip-6pt}
$$
\eqalign{
& g(x,z)\,a\wdg b\wdg y-g(x,y)\,a\wdg b\wdg z\cr
&\qquad=g(x,a)\,b\wdg y\wdg z-g(x,b)\,a\wdg y\wdg z=-Tx\wdg y\wdg z
}
$$
since $\bigwedge^4V=0$ and $-x\wdg Ty\wdg z-x\wdg y\wdg Tz=Tx\wdg y\wdg z$ 
since $\tr T=0$. This yields
$$
Y_1Y_2(tx)-\si Y_2Y_1(xt)=2\,Tx\wdg y\wdg z,
$$
completing the proof.
\endproof

It is immediately clear from Lemma 4.1 that (4.4) is satisfied for 
$t=x\wdg Tx$. This equality is satisfied also for $t=a\wdg b$ since $Tx$ lies 
in the subspace of $V$ spanned by $a$ and $b$. Since 
$Y(x^2)=g(x,x)\,a\wdg b+2\,x\wdg Tx$, equality (4.4) holds for $t=Y(x^2)$ as 
well, and we are done.

\section
5. Final results

We assume that $V$ is a vector space of dimension 3 over an arbitrary field 
$\bbk$ of characteristic $\ne2$. Denote by $\HeckeSym(V)$ the set of all Hecke 
symmetries on $V$ and by $\HeckeSym_0(V)$ its subset consisting of those Hecke 
symmetries for which $\bbS(V,R)=\bbS(V)$, the ordinary symmetric algebra of $V$.

\proclaim
Theorem 5.1.
Each Hecke symmetry $R\in\HeckeSym_0(V)$ with parameter $q$ of the Hecke 
relation is given by the formula
$$
R(xy)={q-1\over2}xy+{q+1\over2}yx-g(x,y)\,a\wdg b-x\wdg Ty-y\wdg Tx,\qquad 
x,y\in V,
$$
where $a,b\in V$ are two vectors, $g:V\times V\to\bbk$ a symmetric bilinear 
form satisfying\/ $(3.1)$ and $T:V\to V$ the linear operator defined in terms of 
$a,b$ and $g$ by\/ $(3.3)$.
\endproclaim

\Proof.
The $R$-skewsymmetrizer $Y$ is determined in Proposition 3.1, and $R$ is found 
by the formula $R=q\cdot\Id-Y$ according to (2.1).
\endproof

The vector space $V^{\ot2}$ is a $\GL(V)$-module in a natural way. Hence there 
is also an action of $\GL(V)$ on the algebra $\End_\bbk V^{\ot2}$. It is clear 
that the set $\HeckeSym(V)$ and its subset $\HeckeSym_0(V)$ are stable under 
this action of $\GL(V)$.

\setitemsize(i4)
\proclaim
Theorem 5.2.
There is a\/ $\GL(V)$-equivariant bijection between the set of Hecke 
symmetries\/ $\HeckeSym_0(V)$ and the set consisting of all pairs $(q,F)$ 
where $0\ne q\in\bbk$ and $F:V^{\ot2}\to V^{\ot2}$ is a linear operator 
satisfying the following conditions:

\item(i1)
\quad
$F(xy)=F(yx)\,$ for all $x,y\in V,$

\item(i2)
\quad
$\Im F\sbs\Alt_2\,,$

\item(i3)
\quad
$\rk F\le1,$

\item(i4)
\quad
$(q-1)^2=-4\mskip1mu\De(F)$

\noindent
where $\De(F)$ is defined by the formula $\,\De(F)=g(a,a)\,g(b,b)-g(a,b)^2\,$ 
with $a,b\in V$ and  a symmetric bilinear form $g:V\times V\to\bbk$ taken so that
$$
F(xy)=g(x,y)\,a\wdg b\qquad\hbox{for all $\,x,y\in V$}.\eqno(5.1)
$$
\endproclaim

\Proof.
Each Hecke symmetry $R\in\HeckeSym_0(V)$ is determined by a symmetric bilinear 
form $g:V\times V\to\bbk$ and a bivector $t=a\wedge b\in\bigwedge^2V$. On the 
other hand, the pair $(g,t)$ gives rise to a linear operator $F$ defined by 
formula (5.1). Clearly $F$ satisfies (i1)--(i3). We can describe a relationship 
between $R$ and $F$ in invariant terms. Let $Y$ be the $R$-skewsymmetrizer (2.1). 
Then
$$
x\wdg Y(xz)=g(x,x)\,a\wdg b\wdg z\qquad\hbox{for all $x,z\in V$}
$$
since the bijection $\omt:\Alt_3\to\bbk$ transforms the above equality into 
the equality
$$
\ell_{xx}(z)=g(x,x)\,\om(a,b,z)
$$
which holds in view of (3.5) and (3.6). Hence $\,F(xx)\wdg z=x\wdg Y(xz)$, and 
linearizing this identity we get
$$
2F(xy)\wdg z=x\wdg Y(yz)+y\wdg Y(xz),\qquad x,y,z\in V.\eqno(5.2)
$$
It follows that $F$ is uniquely determined by $R$. Since the multiplication 
$\wdg$ is $\GL(V)$-invariant, the map given by the assignment $R\mapsto F$ is
$\GL(V)$-equivariant.

We put into correspondence to $R$ the pair $(q,F)$ where $q$ is the parameter 
of the Hecke relation satisfied by $R$ and $F$ is found from (5.2). Then (i4) 
is a consequence of (3.1). By means of the bijection $\omt$ formula (3.7) 
translates into the identity
$$
x\wdg Y(yz)=F(xy)\wdg z+F(xz)\wdg y-F(yz)\wdg x+{q+1\over2}\,x\wdg y\wdg z.
$$
Hence $Y$, and therefore also $R$, are uniquely determined by the pair $(q,F)$.
\endproof

\proclaim
Corollary 5.3.
Given a Hecke symmetry $R\in\HeckeSym_0(V)$ with parameter $q,$ put
$$
R_\la=R_0+\la(R-R_0)
$$
where $R_0:V^{\ot2}\to V^{\ot2}$ is the flip of tensorands $xy\mapsto yx$. 
Then $R_\la\in\HeckeSym_0(V)$ for all $\la\in\bbk$ such that $\,\la(q-1)\ne-1$. 
Thus $R$ is a deformation of $R_0$.
\endproclaim

\Proof.
We use the bijective correspondence described in Theorem 5.2. Let $(q,F)$ be 
the pair corresponding to $R$. Then $R_\la$ is the Hecke symmetry corresponding 
to the pair $(q_\la,\la F)$ where $q_\la=1+\la(q-1)$. Note that
$$
(q_\la-1)^2=\la^2(q-1)^2=-4\mskip1mu\De(\la F)
$$
since $\De(\la F)=\la^2\De(F)$. If $\la$ is such that $\la(q-1)=-1$, then 
$q_\la=0$, and the corresponding operator $R_\la$ is singular. This value of 
$\la$ should be excluded.
\endproof

\proclaim
Corollary 5.4.
If $\,R\in\HeckeSym_0(V),$ then the operator $r=R_0R-\Id$, regarded as an 
element of $\,\gl(V)\ot\gl(V),$ is a solution to the classical Yang-Baxter 
equation. Moreover,
$$
r+r_{21}=(q-1)(R_0+\Id)
$$
where $\,r_{21}=R_0{\mskip1.6mu}rR_0=RR_0-\Id$. Thus $r$ is skewsymmetric 
precisely when $q=1$.
\endproclaim

\Proof.
Since $R_\la$ in Corollary 5.3 satisfies the braid equation, the operator 
$$
R_0R_\la=\Id+\la r
$$
satisfies the constant quantum Yang-Baxter equation, for all $\la\in\bbk$. This 
entails the classical Yang-Baxter equation for $r$ (see \cite{Ch-P}). We have
$$
R_\la^2=\Id+(\la-\la^2)(r+r_{21})+\la^2(R^2-\Id).
$$
Note that $\,R^2-\Id=(q-1)(R+\Id)$. Hence
$$
R_\la^2-\la(q-1)(R_\la+\Id)-\Id=(\la-\la^2)\bigl(r+r_{21}-(q-1)(R_0+\Id)\bigr),
$$
which is 0 since $R_\la$ satisfies the Hecke relation with parameter 
$\,q_\la=1+\la(q-1)$. This proves the desired formula for $\,r+r_{21}$.
\endproof

\section
6. Determination up to equivalence

Suppose that the base field $\bbk$ is algebraically closed of characteristic 
$\ne2$. There is a choice of a basis $x_1,x_2,x_3$ for the vector space $V$ 
which brings any Hecke symmetry $R$ described in Theorem 5.1 to one of 8 
types. If $a\wedge b=0$, then $R=R_0$ is the flip of tensorands. Suppose that 
$a\wedge b\ne0$. The first two basis vectors $x_1,x_2$ will be chosen so that 
$x_1\wedge x_2=a\wedge b$. Changing $a$ and $b$, we may assume that $x_1=a$ 
and $x_2=b$. Clearly, the rank of the bilinear form $g$ and the rank of its 
restriction to the 2-dimensional subspace $\<a,b\>$ are $\GL(V)$-invariants of 
$R$. It is seen from (3.1) that $q=1$ if and only if $g|_{\<a,b\>}$ is 
degenerate.

If $g|_{\<a,b\>}$ is nondegenerate, then $\<a,b\>$ has a basis consisting of 
isotropic vectors. So we may assume in this case that $g(a,a)=g(b,b)=0$, 
and then $4g(a,b)^2=(q-1)^2$ by (3.1). Since $g(-b,a)=-g(a,b)$, we can achieve 
$g(a,b)=(q-1)/2$ by replacing the pair $a,b$ with $-b,a$, if necessary. The 
third basis vector $x_3$ can be chosen in the orthogonal complement of 
$\<a,b\>$. We thus obtain two types of Hecke symmetries with $q\ne1$ 
distinguished by whether $g$ is nondegenerate or not.

If $g|_{\<a,b\>}$ is degenerate, then it either has rank 1 or is identically 
zero, while
$$
\rk g|_{\<a,b\>}\le\rk g\le2+\rk g|_{\<a,b\>}.
$$
If $\rk g|_{\<a,b\>}=1$, then $x_2$ can be taken in the radical of $g|_{\<a,b\>}$, 
and $x_3$ orthogonal to $x_1$. If $g|_{\<a,b\>}=0$, then $x_1$ can be taken in 
the radical of $g$. So $g(x_1,x_3)=0$ in both cases. If $g(x_2,x_3)\ne0$, then 
$x_3$ can be chosen isotropic. Altogether there are 6 types of Hecke 
symmetries with $q=1$.

The matrices of the bilinear forms $g$ in suitably chosen bases of $V$ are 
listed below for all 8 types of Hecke symmetries $R$ with the $R$-symmetric 
algebra $\,\bbS(V)$:
$$
\hbox{Types 1, 2\quad($q\ne1$):}\qquad
M(g)=\pmatrix{0&\displaystyle{q-1\over2}&0\cr\displaystyle{q-1\over2}&0&0\cr0\vphantom{\Bigl|}&0&1}\quad
\pmatrix{0&\displaystyle{q-1\over2}&0\cr\displaystyle{q-1\over2}&0&0\cr0\vphantom{\Bigl|}&0&0}
\vadjust{\vskip-15pt}
$$
Types 3--8\quad($q=1$):
$$
\pmatrix{1&0&0\cr0&0&1\cr0&1&0}\
\pmatrix{1&0&0\cr0&0&0\cr0&0&1}\
\pmatrix{1&0&0\cr0&0&0\cr0&0&0}\
\pmatrix{0&0&0\cr0&0&1\cr0&1&0}\
\pmatrix{0&0&0\cr0&0&0\cr0&0&1}\
\pmatrix{0&0&0\cr0&0&0\cr0&0&0}
$$

It is seen from the above description that equivalence classes of Hecke 
symmetries in the set $\HeckeSym_0(V)$ are distinguished by the parameter $q$ 
and the ranks of $g$ and $g|_{\<a,b\>}$.

The values of $R$ at the basis elements $x_ix_j$ of the space $V^{\ot2}$ can 
be computed readily by means of the explicit formula in Theorem 5.1. In Type 1 
$\,x_1,x_2,x_3$ are eigenvectors of $T$ with respective eigenvalues $(q-1)/2$, 
$-(q-1)/2$, 0, and we get the formulas
$$
\hfuzz=1.8em
\vcenter{\openup2\jot
\halign{\hfil$#$&${}#$\hfil&
\qquad\hfil$#$&${}#$\hfil&\qquad\hfil$#$&${}#$\hfil\cr
R(x_1^2)&=qx_1^2&R(x_1x_2)&=(q-1)x_1x_2+x_2x_1&R(x_1x_3)&=(q-1)x_1x_3+x_3x_1\cr
R(x_2x_1)&=qx_1x_2&R(x_2^2)&=qx_2^2&R(x_2x_3)&=qx_3x_2\cr
R(x_3x_1)&=qx_1x_3&R(x_3x_2)&=(q-1)x_3x_2+x_2x_3&R(x_3^2)&=qx_3^2-x_1x_2+x_2x_1\cr
}}
$$

Let $\{E_{ij}\mid1\le i,j\le3\}$ be the basis of $\gl(V)$ consisting of the 
matrix units with respect to the chosen basis of $V$. Then the classical 
$r$-matrix $r=R_0R-\Id$ described in Corollary 5.4 is written as
$$
\hfuzz=2.4em
\eqalign{
r&=(q-1)(E_{11}\ot E_{11}+E_{22}\ot E_{22}+E_{33}\ot E_{33}
+E_{22}\ot E_{11}+E_{22}\ot E_{33}+E_{33}\ot E_{11})\cr 
&\qquad+(q-1)(E_{21}\ot E_{12}+E_{23}\ot E_{32}+E_{31}\ot E_{13})
+E_{13}\ot E_{23}-E_{23}\ot E_{13}\,.
}
$$

In Type 2 the operator $R$ has the same values at the monomials $x_ix_j$ except 
that $R(x_3^2)=qx_3^2$, and the $r$-matrix is as displayed above, but without 
the last two summands $E_{13}\ot E_{23}$ and $E_{23}\ot E_{13}$. Specializing 
the parameter $q$ to 1 gives the formulas for Types 7 and 8.

In Type 3 we have $\,Tx_1=-x_2$, $Tx_2=0$, $Tx_3=x_1$, which gives
$$
\hfuzz=1.8em
\vcenter{\openup2\jot
\halign{\hfil$#$&${}#$\hfil&
\qquad\hfil$#$&${}#$\hfil&\qquad$#$\hfil\cr
R(x_1^2)&=x_1^2+x_1x_2-x_2x_1&R(x_1x_2)&=x_2x_1&R(x_1x_3)=x_3x_1-x_2x_3+x_3x_2\cr
R(x_2x_1)&=x_1x_2&R(x_2^2)&=x_2^2&R(x_2x_3)=x_3x_2\cr
R(x_3x_1)&=x_1x_3-x_2x_3+x_3x_2&R(x_3x_2)&=x_2x_3&R(x_3^2)=x_3^2+2(x_1x_3-x_3x_1)\cr
\noalign{\medskip}
\noalign{\centerline{and\qquad
$r=E_{21}\wdg(E_{11}+E_{33})+E_{23}\wdg E_{31}+2{\mskip1mu}E_{33}\wdg 
E_{13}\,$.}}
}}
$$

The formulas for Types 4, 5, 6 are almost the same as for Type 3 with only a 
few changes. In Type 4 the differences are in the formulas $Tx_3=0$,
$$
R(x_3^2)=x_3^2-x_1x_2+x_2x_1\quad\hbox{and}\quad
r=E_{21}\wdg(E_{11}+E_{33})+E_{23}\wdg(E_{31}-E_{13}).
$$
In Type 5 we have $\,Tx_3=0$,
$$
R(x_3^2)=x_3^2\quad\hbox{and}\quad
r=E_{21}\wdg(E_{11}+E_{33})+E_{23}\wdg E_{31},
$$
and in Type 6 we have $\,Tx_1=0$,
$$
R(x_1^2)=x_1^2,\quad R(x_1x_3)=x_3x_1,\quad R(x_3x_1)=x_1x_3,\quad 
r=2{\mskip1mu}E_{33}\wdg E_{13}\,.
$$

If $\chr\bbk\ne3$, then the Lie algebra $\gl(3)$ is the direct sum of $\sl(3)$ 
and the onedimensional center. In this case the projection $\gl(3)\to\sl(3)$ 
transforms each classical $r$-matrix in $\gl(3)\ot\gl(3)$ to one in 
$\sl(3)\ot\sl(3)$. All classical $r$-matrices in $\sl(3)\ot\sl(3)$ were 
described by Gerstenhaber and Giaquinto, and for each a quantization was found 
\cite{Ger-Gi98}. Certainly, a large part of those $r$-matrices do not 
correspond to Hecke symmetries with $\bbS(V,R)=\bbS(V)$.

Recall that a finite-dimensional Lie algebra $L$ is said to be 
\emph{quasi-Frobenius} if there exists a nondegenerate alternating bilinear 
form on $L$ which is a 2-cocycle, and if such a form is a 2-coboundary, then 
$L$ is \emph{Frobenius}. It is a general fact that the carriers of 
skewsymmetric classical $r$-matrices are quasi-Frobenius Lie algebras (see 
\cite{Dr83} and \cite{Ch-P, Prop. 2.2.6}). In our present study the carriers 
of the $r$-matrices corresponding to Hecke symmetries of Types 3-6 are the 
following four Frobenius subalgebras of the Lie algebra $\gl(3)$:
$$
\displaylines{
\<E_{11},\,E_{13},\,E_{21},\,E_{23},\,E_{31},\,E_{33}\>,\qquad
\<E_{11}+E_{33},\,E_{13}-E_{31},\,E_{21},\,E_{23}\>,\cr
\<E_{11}+E_{33},\,E_{21},\,E_{23},\,E_{31}\>,\qquad\<E_{13},\,E_{33}\>.
}
$$
The carrier in Type 7 is the abelian subalgebra $\<E_{13},E_{23}\>$, and in 
Type 8 the carrier is 0. One can see that these six Lie algebras are pairwise 
nonisomorphic. The conjugacy classes of all quasi-Frobenius subalgebras in the 
Lie algebra $\sl(3)$ were determined by Stolin \cite{Sto91}.

\references
\nextref Ch-P
\auth{V.,Chari;A.,Pressley}
\book{A Guide to Quantum Groups}
\publisher{Cambridge Univ. Press}
\Year{1994}

\nextref Dr83
\auth{V.G.,Drinfeld}
\paper{Hamiltonian structures on Lie groups, Lie bialgebras and the geometric meaning of the classical Yang-Baxter equations\inRus}
\journal{Dokl. Akad. Nauk SSSR}
\Vol{268}
\Year{1983}
\Pages{285-287};
\etransl{Sov. Math. Dokl.}
\Vol{27}
\Year{1983}
\Pages{68-71}

\nextref Ew-O94
\auth{H.,Ewen;O.,Ogievetsky}
\paper{Classification of the $GL(3)$ quantum matrix groups}\hfill\break
\hbox{}\hfill arXiv:9412009.

\nextref Ger-Gi98
\auth{M.,Gerstenhaber;A.,Giaquinto}
\paper{Boundary solutions of the quantum Yang-Baxter equation and solutions in three dimensions}
\journal{Lett. Math. Phys.}
\Vol{44}
\Year{1998}
\Pages{131-141}

\nextref Gur90
\auth{D.I.,Gurevich}
\paper{Algebraic aspects of the quantum Yang-Baxter equation\inRus}
\journal{Algebra i Analiz}
\Vol{2} (4)
\Year{1990}
\Pages{119-148};
\etransl{Leningrad Math. J.}
\Vol{2}
\Year{1991}
\Pages{801-828}

\nextref Sto91
\auth{A.,Stolin}
\paper{Constant solutions of Yang-Baxter equation for\/ $\sl(2)$ and\/ $\sl(3)$}
\journal{Math. Scand.}
\Vol{69}
\Year{1991}
\Pages{81-88}

\endreferences

\bye